\documentclass [12pt] {article}
\title {A refined 1-cocycle for regular isotopies and the refined tangle equations }
\author{Thomas Fiedler}
\usepackage{amsmath, amsfonts, amssymb, latexsym, graphicx}
\begin{document}

\newcommand{\Tr}{\mathrm{Tr}}
\renewcommand{\L}{\mathbb{L}}
\newtheorem{proposition}{Proposition}
\newtheorem{theorem}{Theorem}
\newtheorem{lemma}{Lemma}
\newtheorem{corollary}{Corollary}
\newtheorem{example}{Example}
\newtheorem{remark}{Remark}
\newtheorem{definition}{Definition}
\newtheorem{question}{Question}
\newtheorem{conjecture}{Conjecture}
\newtheorem{observation}{Observation}
\maketitle  
\begin{center}

\end{center}

\begin{abstract}

We refine the combinatorial  1-cocycle $\mathbb{L}R_{reg}$ for regular isotopies of long knots to a 1-cocycle with values in the free $\mathbb{Z}[x,x^{-1}]$-module generated by regular isotopy classes of oriented tangles with exactly one signed ordinary double point. 

We use it to define the refined tangle equations for couples of knot diagrams, where the coefficients are now Laurent polynomials instead of integers.  A solution of the tangle equations gives quantitative information about any knot isotopy which relates two given knot diagrams.  If the tangle equations have no solution, then the diagrams represent different knots.

\end{abstract}

\tableofcontents

\section{Introduction and main results}

The {\em combinatorial  1-cocycle} $\mathbb{L}R_{reg}$ for regular isotopies of long knots and their 2-cables was introduced in \cite{FMVZ}.  This 1-cocycle takes its values in the free $\mathbb{Z}$-module generated by regular isotopy classes of oriented tangles with exactly one signed ordinary double point.  Its evaluation on the {\em push arc} in the moduli space (i.e.  pushing a small auxiliary long knot diagram from one end of a long knot to the other) has lead to the {\em tangle equations}.  They were generalized for knotoids in \cite{FKL}.  A weaker form of the 1-cocycle for long knots and their cables, where the double point is replaced by its orientation preserving smoothing and the sign of the double point is ignored,  was already contained in \cite{F} and has given the {\em quantum equations} in \cite{FZ}.

The 1-cocycle $\mathbb{L}R_{reg}(x)$ can be seen as an {\em elementary quantization} of the 1-cocycle $\mathbb{L}R_{reg}$.

For the convenience of the reader we have kept the present paper essentially self-contained and we have tried to keep it relatively short.
 \vspace{0,4cm}

We fix once for all a linear projection $pr: (x,y,z) \to (x,y)$ of the oriented 3-space into the oriented 2-plan. 

Let $D$ be a diagram of an oriented long knot.   If we smooth a crossing with respect to the orientation of $D$, then we obtain an oriented circle and a long knot.  There are two global types of crossings: for a {\em 0-crossing} the circle goes from its over-cross to its under-cross and for a {\em 1-crossing} it goes from the under-cross to its over-cross.  We denote the sum of the usual signs of the $1$-crossings of $D$ by $w_1(D)$ and we denote the sum of the signs of all crossings of $D$ by $w(D)$.  This are finite type knot invariants of degree 1 for regular isotopy (for finite type knot invariants and finite type 1-cocycles see \cite{V1} and \cite{V2}).

Let $M$ be the (disconnected) topological moduli space of all oriented diagrams of long knots with respect to regular isotopy (i.e. without Reidemeister I moves, for the Reidemeister moves see e.g. \cite{K}) and let $M_0$ be the topological moduli space of these diagrams and which have in addition exactly one closed oriented component.  

Let $M_+$ (respectively $M_-$)  be the topological moduli spaces of all diagrams of oriented long knots with exactly one positive signed (respectively negative signed) ordinary double point in the plan and which stays in all regular  isotopies an ordinary double point in the plan (in other words a vertex rigid isotopy of the double point). 
 \vspace{0,4cm}

Let $D$ and $D'$ be two long knot diagrams with the same Whitney index and the same writhe $w(D)=w(D')$ (see e.g.  \cite{K}) and let $K$ be an auxiliary long knot diagram with given $w_1(K)$. The following proposition is a preparation for our main result.

\begin{proposition}

There is a combinatorial 1-cocycle 

$\mathbb{L}R_{reg}(x) \in Z^1(M; H_0(M_+;\mathbb{Z}[x,x^{-1}]) \bigoplus  H_0(M_-;\mathbb{Z}[x,x^{-1}]))$.  
 \vspace{0,2cm}

Moreover,  if $D$ and $D'$ are regularly isotopic, then
 \vspace{0,2cm}

$\mathbb{L}R_{reg}(x)(push(K,D)) - \mathbb{L}R_{reg}(x)(push(K,D'))=$

 $K \sum_i D_i(\sum_j  a_{i,j}(x^{b_{i,j}+w_1(K)}-x^{b_{i,j}}))$.
 \vspace{0,2cm}

Here $a_{i,j}$ and $b_{i,j}$ are integers and $D_i$ are regular isotopy classes of singular long knots with exactly one double point.  This is the refined tangle equation, which moreover splits with respect to the positive and the negative double points.

\end{proposition}

$\mathbb{L}R_{reg}(x)$ will be defined in an explicit way in the next section.

 \vspace{0,4cm}

We  double now the long knot $D$ by using the blackboard framing and push the 2-cable of $K$, denoted by $2K$, along the 2-cable of $D$, denoted by $2D$.  We distinguish the knot drawn in red, from its parallel longitude, drawn in black.  

Analogous to the previous definitions, let $M^{red}$ be the topological moduli space with respect to regular isotopy of all red oriented long knot diagrams $D$ together with their parallel black longitude and let $M^{red}_0$ be the topological moduli space of couples of long knots and which have in addition exactly one closed oriented component.  Let $M^{red}_+$ (respectively $M^{red}_-$) be the topological moduli space of all oriented red singular long knot diagrams $D$ together with their black parallel longitude and with exactly one positive signed (respectively negative signed) ordinary double point in the plan and such that both branches  in the double point come from the red component $D$.  The double points in the image of the 1-cocycles come always from some $0$-crossings of the long knots and the sign of the double point is just the sign of the corresponding $0$-crossing.

Completely analogous we define $M^{black-red}_+$ and $M^{black-red}_-$ as the topological moduli space of all oriented red singular long knot diagrams $D$ together with their black parallel longitude and with exactly one positive signed (respectively negative signed) ordinary double point in the plan and such that the under-crossing of the corresponding $0$-crossing for the double point is on the black longitude and the over-crossing is on the red knot, see the next section.

Let $K$ be an auxiliary long knot diagram with given $w_1(K)+w(K)$.

The following theorem is our main result.

\begin{theorem}

There are combinatorial 1-cocycles 

$\mathbb{L}R^{red}_{reg}(x) \in Z^1(M^{red}; H_0(M^{red}_+;\mathbb{Z}[x,x^{-1}]) \bigoplus  H_0(M^{red}_-;\mathbb{Z}[x,x^{-1}]))$ and 

$\mathbb{L}R^{black-red}_{reg}(x) \in Z^1(M^{red}; H_0(M^{black-red}_+;\mathbb{Z}[x,x^{-1}]) \bigoplus  H_0(M^{black-red}_-;\mathbb{Z}[x,x^{-1}]))$.
 \vspace{0,2cm}

Moreover,  if $D$ and $D'$ are regularly isotopic, then
 \vspace{0,2cm}

$\mathbb{L}R^{red}_{reg}(x)(push(2K,2D)) - \mathbb{L}R^{red}_{reg}(x)(push(2K,2D'))=$

 $2K \sum_i D_i(\sum_j  a_{i,j}(x^{b_{i,j}+w_1(K)+w(K)}-x^{b_{i,j}}))$ and 
 \vspace{0,2cm}

$\mathbb{L}R^{black-red}_{reg}(x)(push(2K,2D)) - \mathbb{L}R^{black-red}_{reg}(x)(push(2K,2D'))=$

 $2K \sum_i D_i(\sum_j  a_{i,j}(x^{b_{i,j}+w_1(K)+w(K)}-x^{b_{i,j}}))$. 
 \vspace{0,2cm}

Here $a_{i,j}$ and $b_{i,j}$ are again integers (which depend of course of the type red or black-red of the 1-cocycle) and $D_i$ are here regular isotopy classes of oriented singular long 2-tangles with exactly one double point which comes from a red-red or correspondingly black-red crossing.  This are the refined tangle equations for long knots with their longitude,  and which moreover split with respect to the positive and the negative double points.

\end{theorem}

If the refined tangle equations have no solution, then $D$ and $D'$ represent different knots.  Notice that the numbers $a_{i,j}$ and $b_{i,j}$ and the regular isotopy classes of the $D_i$ do not depend on the chosen auxiliary long knot diagrams $K$.

In the degenerate case of $w_1(K)+w(K)=0$, it follows from the theorem that $\mathbb{L}R^{red}_{reg}(x)(push(2K,2D))$ and $\mathbb{L}R^{black-red}_{reg}(x)(push(2K,2D))$ are  knot invariants of $D$ and which depend on $K$ as a parameter. (As well known,  two oriented classical knots are isotopic in $S^3$ if and only if the corresponding long knots share the same Whitney index and the same writhe and are regularly isotopic as long knots.) We explain in Section 5 in much detail why we think that this is a non-trivial knot invariant!

$\mathbb{L}R^{red}_{reg}(x)(push(2K,2D))$ and $\mathbb{L}R^{black-red}_{reg}(x)(push(2K,2D))$ are linear combinations with polynomial coefficients in $x$ of 2-cables of  long knots,  each with exactly one signed ordinary double point.  We could apply now the invariant $\rho_1$ for singular knots from \cite{SvV} to them.

Let us define two morphisms in the following way for $red$ and also for $black-red$:
\vspace{0,2cm}

$os: H_0(M_+;\mathbb{Z}[x,x^{-1}]) \bigoplus  H_0(M_-;\mathbb{Z}[x,x^{-1}]) \to H_0(M_0;\mathbb{Z}[x,x^{-1}])$
\vspace{0,2cm}

$cc: H_0(M_+;\mathbb{Z}[x,x^{-1}]) \bigoplus  H_0(M_-;\mathbb{Z}[x,x^{-1}]) \to H_0(M;\mathbb{Z}[x,x^{-1}])$
\vspace{0,2cm}

In $os$ we smooth all double points with respect to the orientation (this gives each time a 2-string link and a closed component).  In $cc$ we replace each double point in $H_0(M_+;\mathbb{Z}[x,x^{-1}])$ by a negative crossing and each double point in $H_0(M_-;\mathbb{Z}[x,x^{-1}])$ by a positive crossing.  This will correspond to changing each time a crossing  to its opposite in the 2-cable diagram $2K2D$.  The important point is, that after applying $os$ or $cc$ we obtain 2-string links (with perhaps a closed component in addition) which are no longer 2-cables of a long knot!

There aren't known any finite type invariants which could detect the orientation of a long knot. However,  Duzhin and Karev have shown in \cite{DK}, that there is a finite type invariant of degree 7 for 2-component string links which can sometimes detect the orientation of the 2-string link if it is not just the 2-cable of a long knot.  We can apply now Duzhin and Karevs invariant for $red$ and $black-red$ to $cc(\mathbb{L}R_{reg}(x)(push(2K,2D)))$ and to the 2-string part of $os(\mathbb{L}R_{reg}(x)(push(2K,2D)))$ as well. This could perhaps already detect the orientation of $D$.  

But perhaps already the tangle equations are enough to distinguish a knot from its inverse , because the weights $W_1(d)$ and $W_1(ml)$ for a move in $2D$ and the corresponding move in its inverse $2D'$ are very different.  Notice that the inversion interchanges the crossings $ml$ and $hm$ in the Reidemeister  III moves.

We need very much a computer program in order to calculate examples.
\vspace{0,4cm}

If $D$ and $D'$ are isotopic knots, then the refined tangle equations give quantitative information about any regular isotopy which connects $D$ with $D'$, namely the singular tangles $D_i$ and the numbers $a_{i,j}$ and $b_{i,j}$.  This follows from the fact that we can calculate the left side of the equation which does not depend on the chosen isotopy $\gamma$ from $D$ to $D'$.

\section{The refined 1-cocycles}

For the convenience of the reader we will repeat here the main definitions, compare also \cite{F},  \cite{FZ}, \cite{FMVZ} and \cite{FKL}.

A {\em Gauss diagram of\/} a long knot $D$ is an oriented circle with oriented chords attached and with a marked point $\infty$,  the point at infinity on the knot. The chords correspond to the crossings of the knot and are always oriented from the under-cross to the over-cross (here we use the orientation of the $z$-axes). Moreover, each chord (or arrow) has a sign, denoted by $w(q)$, which corresponds to the usual sign of the crossing $q$.  A Reidemeister III move corresponds to a triangle in the Gauss diagram.

\begin{figure}
\centering
\includegraphics{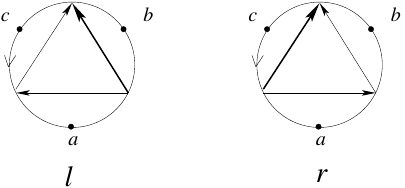}
\caption{\label{globtricross}  The six global types of triple crossings e.g. $r_a$ means the type on the right with $a=\infty$ and $l_c$ means the type on the left with $c=\infty$.}  
\end{figure}

To each Reidemeister move of type III corresponds a diagram with a {\em triple 
crossing} $p$: three branches of the knot (the highest, middle and lowest with respect to the projection $pr$. A small perturbation of the triple crossing leads to an ordinary diagram with three crossings near $pr(p)$.
 We name the corresponding branches of the knot the lowest, middle and highest branch and the corresponding crossings $ml$ (lowest with middle), $hm$ (middle with highest) and $d$ (for distinguished). For better visualization we draw the crossing $d$ always with a thicker arrow. There are exactly six  global types of triple crossings with respect to the point $\infty$, shown in Fig.~\ref{globtricross} . The two new crossings from a Reidemeister II move are also called distinguished and they are denoted by $d$ too.

\begin{figure}
\centering
\includegraphics{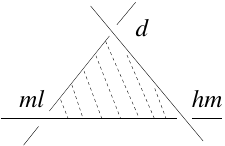}
\caption{\label{names} The names of the crossings in a R III-move}  
\end{figure}

\begin{figure}
\centering
\includegraphics{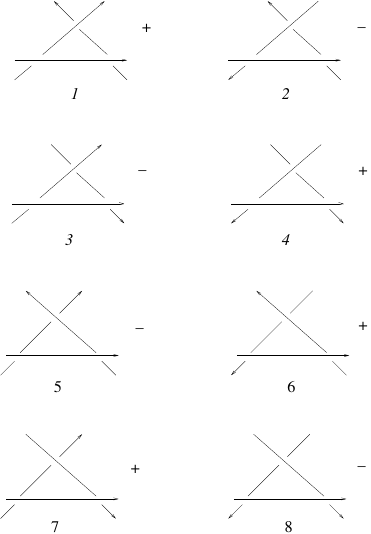}
\caption{\label{loctricross} Local types of a triple crossing together with the local side of the discriminant. The triple crossings of the types 2 and 6 are called {\em star-like}.}  
\end{figure}

There are exactly eight local types of triple crossings. They are shown in Fig.~\ref{loctricross}.  The sign in the figure defines the local side of the discriminant (i.e. the diagrams with exactly one ordinary triple point in the projection). We define it for the local type 1 (i.e.  all three crossings are positive) and the cube equations determine automatically the sides of all other local types, see  \cite{F} or \cite{FMVZ}.

\begin{definition}

The sign of a R III move $p$, denoted by $sign(p)$,  is $+1$ if the move goes from the negative side to the positive side and $-1$ otherwise. 

The sign of a R II move $p$ is $+1$ if it increases the number of crossings and $-1$ otherwise. 

\end{definition}

Let $p$ be a R III or R II move of $D$ where the distinguished crossing $d$ is a $0$-crossing. The Gauss diagram of $D$ with the triple crossing or auto-tangency from $p$ is denoted by $D(p)$.

\begin{definition}
A crossing $q$ of $D(p)$ is called a {\em f-crossing for $d$} if $q$ is of type 1 and  the under-cross of $q$ is in the open oriented arc from $\infty$ to the over cross of $d$ in $D(p)$.

Let $p$ be a triple crossing of type $r_a$ or $l_c$.  A crossing $q$ of $D(p)$ is called a {\em f-crossing for $ml$} if $q$ is of type 1 and  the under-cross of $q$ is in the open oriented arc from $\infty$ to the over cross of $ml$ in $D(p)$.
\end{definition}

\begin{figure}
\centering
\includegraphics{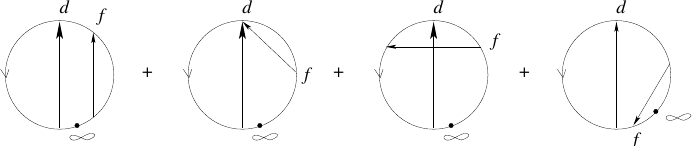}
\caption{\label{foot}  the f-crossings for $d$}  
\end{figure}

\begin{figure}
\centering
\includegraphics{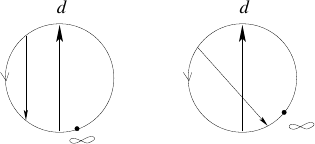}
\caption{\label{notfoot} crossings of type 1 which are not f-crossings for $d$}  
\end{figure}

In other words, the "foot" of $q$ in the Gauss diagram is in the sub arc of the circle which goes from $\infty$ to the head of $d$ and the head of $q$ is in the sub arc which goes from the foot of $q$ to $\infty$.   We illustrate this in Fig.~\ref{foot} and Fig.~\ref{notfoot}.  Notice that we make essential use of the fact that no crossing can ever move over the point at infinity.

\begin{definition}
Let $p$ be a triple crossing of one of the types $r_a$, $r_b$, $l_b$  or let $p$ be an auto-tangency where the distinguished crossing $d$ is a $0$-crossing.  The linear weight $W_1(d)$ is defined as the sum of the signs of all $f$-crossings for $d$.

Let $p$ be a triple crossing of type $r_a$ or $l_c$.  The linear weight $W_1(ml)$ is defined as the sum of the signs of all $f$-crossings for $ml$.
\end{definition}

\begin{figure}
\centering
\includegraphics{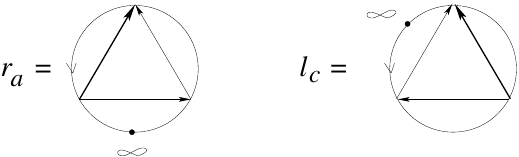}
\caption{\label{fglob} the two global types of triple crossings for which the crossing $ml$ contributes too}  
\end{figure}

The crossing $d$ is a 0-crossing exactly for the types $r_a$,$r_b$ and $l_b$. The crossing $ml$ is a 0-crossing exactly for the types $r_a$, $l_c$ and $l_b$.  But notice that $ml$ does never contribute to the 1-cocycles for the type $l_b$.
 \vspace{0,4cm}

We define now the {\em partial singularizations} of the diagram for the Reidemeister moves. Each singularization preserves the orientations.  We denote a positive double point (it comes from a positive crossing) by $t^+$ and a negative double point by $t^-$.  If the double point is locally in a braid,  then we attach to it its index $i$ as for the braid generators.  

{\em All the singularizations have their values in the corresponding module $H_0(M_+)$ respectively $H_0(M_-)$ .}
\vspace{0,2cm}

{\em Reidemeister II moves }:  Let $p$ be a Reidemeister II move, no matter with equal or opposite tangent direction. If the distinguished crossing $d$ is a 0-crossing, then we replace locally the auto-tangency by   
\vspace{0,4cm}

$II_d(p)= (singularization \thinspace of \thinspace the \thinspace  positive \thinspace crossing)$

$ - (singularization \thinspace  of \thinspace  the \thinspace negative \thinspace  crossing)$.
\vspace{0,4cm}

For example,  if the R II move is of the form $\sigma_1 \sigma_1^{-1}$ then we replace it by $t_1^+\sigma_1^{-1} -\sigma_1t_1^-$. 
\vspace{0,4cm}

{\em Reidemeister III moves}: If the distinguished crossing $d$ is a 0-crossing then we replace the triple crossing by
\vspace{0,4cm}

$III_d(p)= w(d)((singularization  \thinspace  of  \thinspace $d$ \thinspace on \thinspace the\thinspace  positive  \thinspace local \thinspace side \thinspace of \thinspace $p$)$

$ - (singularization \thinspace of \thinspace $d$ \thinspace on \thinspace the \thinspace negative  \thinspace local \thinspace side \thinspace of \thinspace $p$))$
\vspace{0,4cm}

\begin{figure}
\centering
\includegraphics{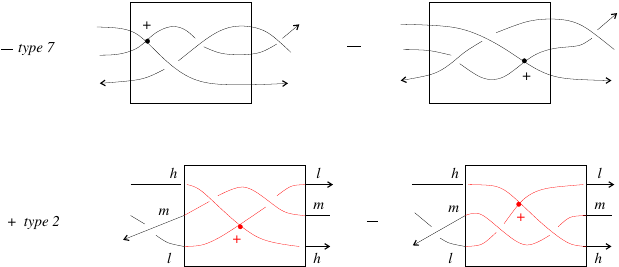}
\caption{\label{eq7-2} The partial $d$-singularizations for the local type 2}   
\end{figure}

\begin{figure}
\centering
\includegraphics{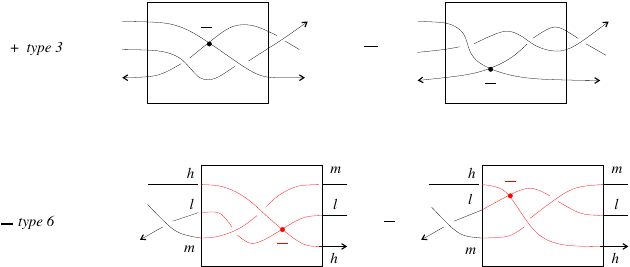}
\caption{\label{eq3-6} The partial $d$-singularizations for the local type 6}  
\end{figure}

For the convenience of the reader we give explicitly $III_d(p)$ for all local cases,  where we write the triple crossing always on the negative side of the discriminant.

type 1 $=\sigma_1\sigma_2\sigma_1$:  $\sigma_2t_1^+\sigma_2-\sigma_1t_2^+\sigma_1$

type 3 $=\sigma_1\sigma_2^{-1}\sigma_1^{-1}$: $\sigma_1\sigma_2^{-1}t_1^--t_2^-\sigma_1^{-1}\sigma_2$

type 4 $=\sigma_1^{-1}\sigma_2^{-1}\sigma_1$: $t_1^-\sigma_2^{-1}\sigma_1-\sigma_2\sigma_1^{-1}t_2^-$

type 5 $=\sigma_2\sigma_1\sigma_2^{-1}$: $\sigma_1^{-1}\sigma_2t_1^+-t_2^+\sigma_1\sigma_2^{-1}$

type 7 $=\sigma_2^{-1}\sigma_1\sigma_2$: $t_1^+\sigma_2\sigma_1^{-1}-\sigma_2^{-1}\sigma_1t_2^+$

type 8 $=\sigma_2^{-1}\sigma_1^{-1}\sigma_2^{-1}$: $\sigma_2^{-1}t_1^-\sigma_2^{-1}-\sigma_1^{-1}t_2^-\sigma_1^{-1}$

star-like type 2 ($w(d)=+1$): the boxes in Fig.~\ref{eq7-2}

star-like type 6 ($w(d)=-1$): the boxes in Fig.~\ref{eq3-6}
\vspace{0,4cm}

If the triple crossing $p$ is of the global type $r_a$ or $l_c$ then we define in addition 
\vspace{0,4cm}

$III_{ml}(p)= (singularization \thinspace of \thinspace $ml$ \thinspace of \thinspace $p$)$
\vspace{0,4cm}

type 1 $=\sigma_1\sigma_2\sigma_1$:  $t_1^+\sigma_2 \sigma_1$

type 3 $=\sigma_1\sigma_2^{-1}\sigma_1^{-1}$: $\sigma_1t_2^-\sigma_1^{-1}$

type 4 $=\sigma_1^{-1}\sigma_2^{-1}\sigma_1$: $\sigma_1^{-1}\sigma_2^{-1}t_1^+$

type 5 $=\sigma_2\sigma_1\sigma_2^{-1}$: $\sigma_2\sigma_1t_2^-$

type 7 $=\sigma_2^{-1}\sigma_1\sigma_2$: $\sigma_2^{-1}t_1^+\sigma_2$

type 8 $=\sigma_2^{-1}\sigma_1^{-1}\sigma_2^{-1}$: $t_2^-\sigma_1^{-1}\sigma_2^{-1}$

star-like type 2 ($w(d)=+1$): the negative crossing $ml$

star-like type 6 ($w(d)=-1$): the positive crossing $ml$
\vspace{0,4cm}

Evidently,  the $ml$-singularization  does not depend on the side here.
\vspace{0,4cm} 

We can now define our tangle-valued 1-cochain.   

\begin{definition} For any generic arc $s \subset M$ the 1-cochain 

$\mathbb{L}R_{reg}(x)(s) \in C^1(M;H_0(M_+;\mathbb{Z}[x,x^{-1}]) \bigoplus  H_0(M_-;\mathbb{Z}[x,x^{-1}]))$ is defined  
\vspace{0,4cm}

$\mathbb{L}R_{reg}(x)(s)= \sum_{p \in s;r_a,l_c} sign(p)(x^{W_1(ml)+(w(hm)-1)/2+1}-x^{W_1(ml)+(w(hm)-1)/2}) III_{ml}(p)$
 \vspace{0,2cm}

$+ \sum_{p \in s; r_a,r_b,l_b}sign(p) x^{W_1(d)}III_d(p)+\sum_{p \in s}sign(p) x^{W_1(d)}II_d(p)$.

\end{definition}

The first part of Proposition 1 says that the cochain $\mathbb{L}R_{reg}(x)$ is a 1-cocycle. This will be proven in the next section.

\begin{remark}
 The 1-cocycle $\mathbb{L}R_{reg}$ from \cite{FMVZ} is just the derivative  $d/dx(\mathbb{L}R_{reg}(x))$ evaluated at $x=1$.  Notice, that the weights $W_1(ml)$ do no longer enter in the definition of $\mathbb{L}R_{reg}$.

The 1-cocycle $R^{(1)}_{reg}$ from \cite{FZ} is just the the 1-cocycle $\mathbb{L}R_{reg}$, where the signs of the double points are ignored, the double points are replaced by their orientation preserving  smoothing,  and each link in the image of the 1-cocycle is replaced by its HOMFLYPT polynomial.

\end{remark}

In order to define our 1-cocycles for the 2-cable we consider the Gauss diagram, where we have glued abstractly to a circle the red long knot with its black long longitude.  There are now two points at infinity in the Gauss diagram,  but we declare the starting point of the red knot as the distinguished point $\infty$.  We define the 1-cocycles exactly as before, but either only with $d$- and $ml$-singularizations of $0$-crossings of  the red knot with itself or alternatively of only with black under-cross and red over-cross (such crossings are automatically $0$-crossings for $\infty$).

The first part of Theorem 1 says that the corresponding cochains $\mathbb{L}R^{red}_{reg}(x)$ and $\mathbb{L}R^{black-red}_{reg}(x)$ are 1-cocycles. We explain in Section 4 how to use them as tangle equations and we indicate in Section 5 why it is likely that they can indeed distinguish knots.

\section{Proof}

The proof is similar to the proof in \cite{FMVZ} that $\mathbb{L}R_{reg}$ is a 1-cocycle , which generalizes the proof in \cite{F} that $R^{(1)}_{reg}$ is a 1-cocycle.

We have to show that $\mathbb{L}R_{reg}(x)$ satisfies the {\em positive tetrahedron equation, the cube equations and the commutation relations}, compare \cite{F}.  The latter are evidently satisfied, because the weights are of degree 1 and hence invariant under R III-moves.  The two new crossings from a R II-move have the same marking, are $f$-crossings only simultaneously but have different signs and cancel out together.

We show that $\mathbb{L}R_{reg}(x)$ satisfies the positive tetrahedron equation in just one case, namely $I_4$,  as an example.  The other case are completely analogous generalizations of the proof for $\mathbb{L}R_{reg}$ and are left to the reader.

The tetrahedron equation is called {\em positive} if all involved six ordinary crossings are positive. Hence, each of the  R III-moves corresponds locally to a braid relation $\sigma_1 \sigma_2 \sigma_1=\sigma_2 \sigma_1 \sigma_2$ in the braid group $B_3$. Moreover, we can attach a sign to the R III-move by declaring that the move has $sign=+1$  if it goes from the expression on the left side to the one on the right side. 

We name the four local branches in the quadruple crossing from the lowest to the highest (with respect to the $z$-coordinate for the projection into the plane) by 1 up to 4, and we name the unique crossing between two branches by the corresponding couple of numbers as well as the triple crossings by the corresponding triple of numbers.

\begin{figure}
\centering
\includegraphics{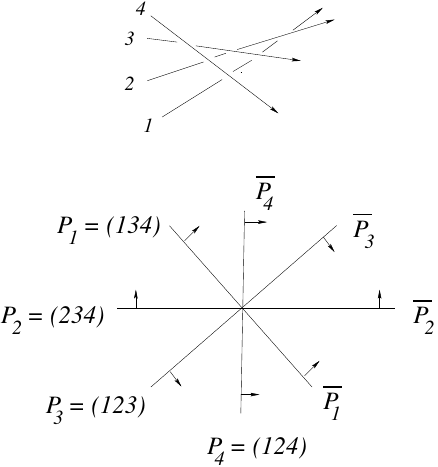}
\caption{\label{quadcros} The names of the strata of triple crossings in the intersection of a normal 2-disc of a positive quadruple crossing}  
\end{figure}

\begin{figure}
\centering
\includegraphics{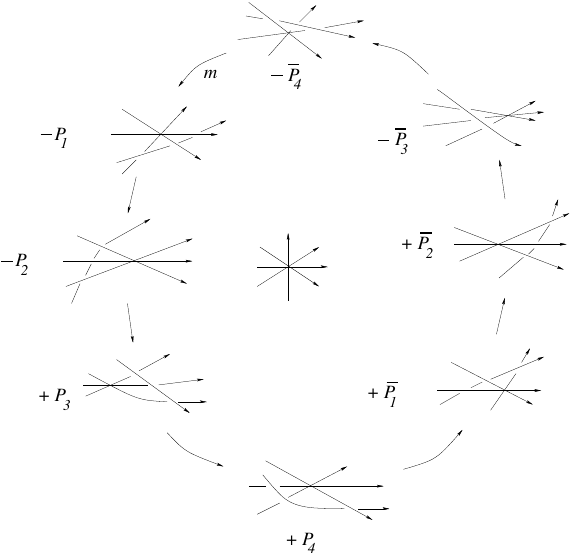}
\caption{\label{unfoldquad} Meridian in the unfolding of a positive quadruple crossing}  
\end{figure}

We show the adjacent strata of positive triple crossings together with their co-orientations (to the positive side) for a positive quadruple crossing in Fig.~\ref{quadcros}. Its unfolding is shown in Fig.~\ref{unfoldquad}, compare  Chapter 3 in \cite{F}.

We denote a positive double point (it comes from a positive crossing) by $t^+$ and a negative double point by $t^-$.  If the double point is locally in a braid,  then we attach to it its index $i$ as for the braid generators.  

\begin{figure}
\centering
\includegraphics{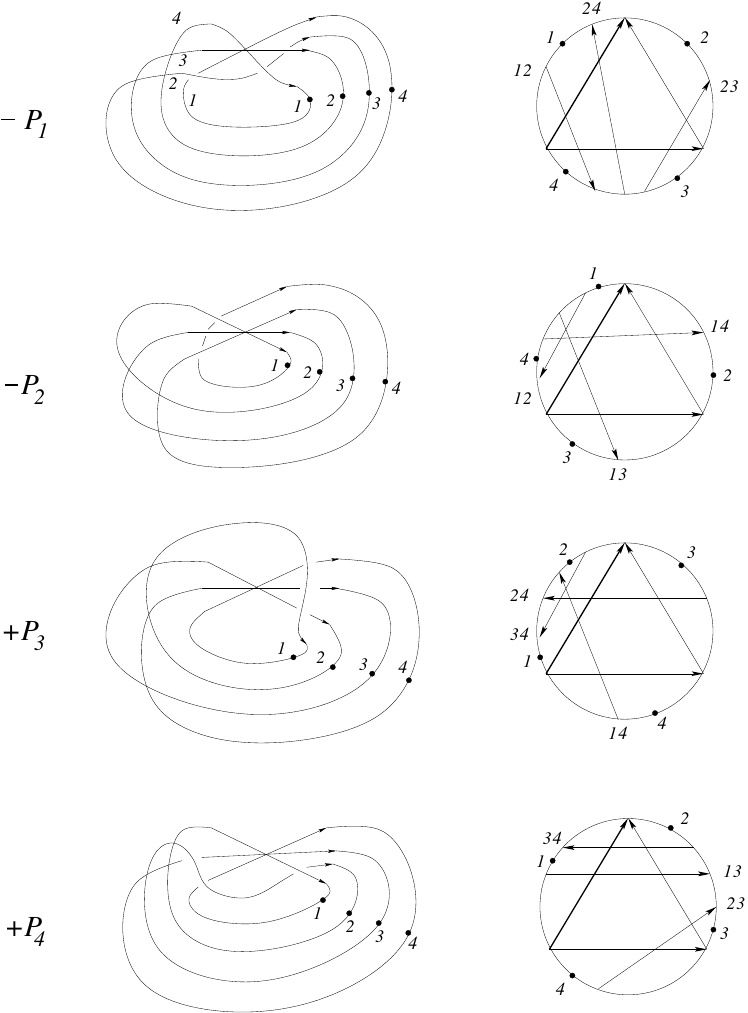}
\caption{\label{Iglob} first half of the meridian for global type I}  
\end{figure}

\begin{figure}
\centering
\includegraphics{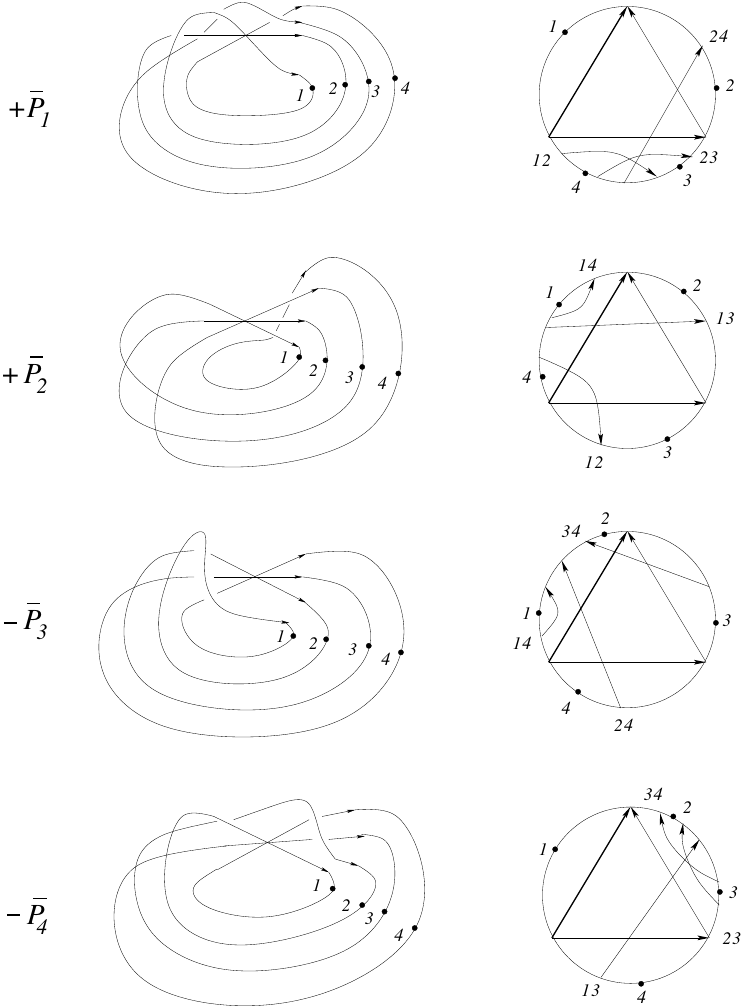}
\caption{\label{I2glob} second half of the meridian for global type I}  
\end{figure}

Fig.~\ref{Iglob} and Fig.~\ref{I2glob} show the eight triple crossings in the meridian of the positive quadruple crossing I (i.e.  all ordinary crossings are positive) where we consider the case $4=\infty$.  In this case $P_1$, $P_3$ and $P_4$ are of type $r_a$ and contribute with $d$- and $ml$-singularizations ($P_2$ is of type $r_c$ and does not contribute at all).  The crossing $34$ is a f-crossing in $-\bar P_3$ but no longer in $P_3$. All the other corresponding weights are just the same and we write $W_1(d)$ and respectively $W_1(ml)$ for the contributions of the invisible crossings in the figure.
\vspace{0,4cm} 

$-P_1$: $-x^{W_1(d=14)+3}(\sigma_1\sigma_3t_2^+\sigma_3\sigma_1\sigma_2-\sigma_1\sigma_2t_3^+\sigma_2\sigma_1\sigma_2)$

$-(x^{W_1(ml=13)+3}-x^{W_1(ml=13)+2})\sigma_1t_2^+\sigma_3\sigma_2\sigma_1\sigma_2$
\vspace{0,4cm}

$+\bar P_1$: $+x^{W_1(d=14)+3}(\sigma_2\sigma_3\sigma_2t_1^+\sigma_2\sigma_3-\sigma_2\sigma_3\sigma_1t_2^+\sigma_1\sigma_3)$

$+(x^{W_1(ml=13)+3}-x^{W_1(ml=13)+2})\sigma_2\sigma_3t_1^+\sigma_2\sigma_1\sigma_3$
\vspace{0,4cm}

$+P_3$: $+x^{W_1(d=13)+2}(\sigma_2t_1^+\sigma_2\sigma_3\sigma_2\sigma_1 - 
\sigma_1t_2^+\sigma_1\sigma_3\sigma_2\sigma_1)$ 

$+(x^{W_1(ml=12)+1}-x^{W_1(ml=12)})t_1^+\sigma_2\sigma_1\sigma_3\sigma_2\sigma_1$
\vspace{0,4cm}

$-\bar P_3$: $-x^{W_1(d=13)+3}(\sigma_3\sigma_2\sigma_1\sigma_3t_2^+\sigma_3 - 
\sigma_3\sigma_2\sigma_1\sigma_2t_3^+\sigma_2)$ 

$-(x^{W_1(ml=12)+2}-x^{W_1(ml=12)+1})\sigma_3\sigma_2\sigma_1t_1^+\sigma_2\sigma_1$
\vspace{0,4cm}

$+P_4$: $+x^{W_1(d=14)+3}(\sigma_2\sigma_1\sigma_3t_2^+\sigma_3\sigma_1-\sigma_2\sigma_1\sigma_2t_3^+\sigma_2\sigma_1)$

$+(x^{W_1(ml=12)+2}-x^{W_1(ml=12)+1})\sigma_2\sigma_1t_2^+\sigma_3\sigma_2\sigma_1$
\vspace{0,4cm}

$-\bar P_4$: $-x^{W_1(d=14)+3}(\sigma_3\sigma_2t_1^+\sigma_2\sigma_3\sigma_2-\sigma_3\sigma_1t_2^+\sigma_1\sigma_3\sigma_2)$

$-(x^{W_1(ml=12)+1}-x^{W_1(ml=12)})\sigma_3t_1^+\sigma_2\sigma_1\sigma_3\sigma_2$
\vspace{0,4cm}

\begin{figure}
\centering
\includegraphics{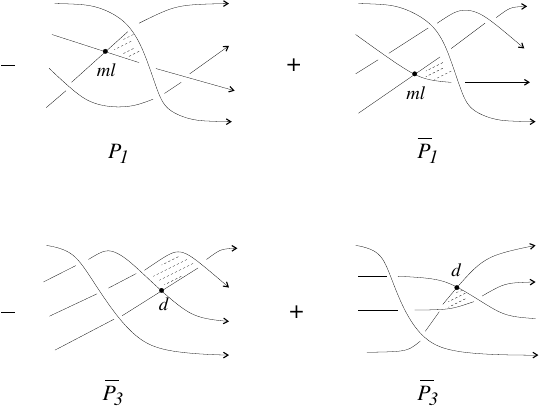}
\caption{\label{P1P3} the $ml$-singularizations in $P_1$ and $\bar P_1$ cancel out with the $d$-singularizations of $P_3$ and $\bar P_3$ because of the coefficients}  
\end{figure}

The $d$-singularizations of $P_1$, $\bar P_1$,  $P_4$ and $\bar P_4$ have all the same coefficient  and cancel out together as for the basic solution $\mathbb{L}$ in \cite{FMVZ}.  The $ml$-singularizations of $P_3$,  $\bar P_3$ and $P_4$, $\bar P_4$ cancel out pairwise,  because the double point from the crossing $12$ is below all other branches and gives no obstruction for isotopies.  The $d$ singularizations from $P_3-\bar P_3$ cancel out with the $ml$-singularizations from $-P_1+\bar P_1$ as shown in Fig.~\ref{P1P3}, for each of the coefficients $x^{W_1(d=ml=13)+3}$ and $x^{W_1(d=ml=13)+2}$.
\vspace{0,4cm}

We have now to solve the cube equations.

\begin{figure}
\centering
\includegraphics{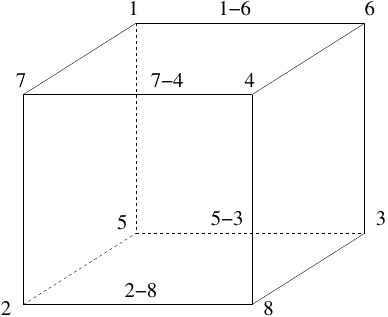}
\caption{\label{gamma} the cube $\Gamma$}  
\end{figure}

The cube is shown in Fig.~\ref{gamma} and the local types of triple crossings were shown in Fig.~\ref{loctricross}.  We have to show that the contributions in the meridian of each edge cancel out together, see \cite{F}.

There are two basic local types of triple crossings for oriented tangles: we call them {\em braid-like} (the only sort which occurs for braids) and {\em star-like} (which never occur for braids),  compare  Fig.~\ref{loctricross}. 

The triple crossings of local types $7,4,8,2$ in the vertical square in front of the cube are exactly the triple crossings for which $w(hm)=-1$.  We will show that $\mathbb{L}R_{reg}(x)$ satisfies the cube equations just for the edge $r1-7$, which explains the correction term $(w(hm)-1)/2$,   for the edge $r7-4$ and for the edge $r1-5$,  where sometimes R II moves contribute non trivially too,  and cancel out either with the $d$- or the $ml$-singularization from the R III moves. The considerations for all other edges, see \cite{F}, are completely analogous and are left to the reader.

\begin{figure}
\centering
\includegraphics{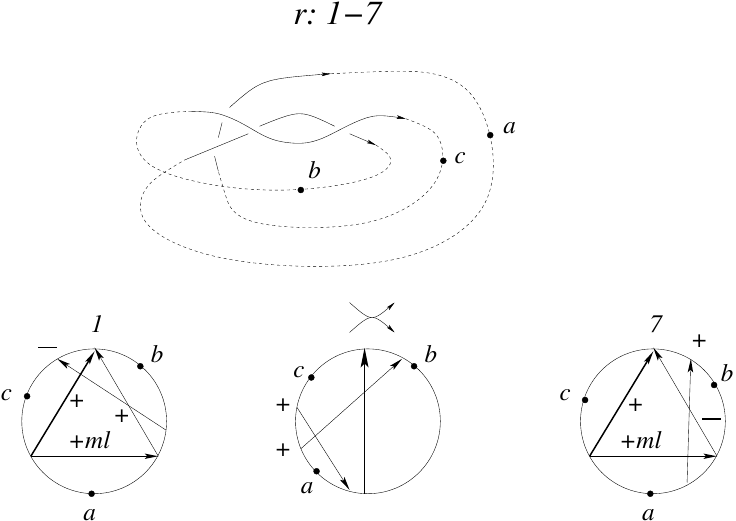}
\caption{\label{r1-7} the edge $r1-7$ in the cube $\Gamma$}  
\end{figure}

The edges $r1-7$ is shown in Fig.~\ref{r1-7}.

For the global type $r_a$ the crossing $hm$ contributes to $W_1(d)$ and hence for all edges of the cube $W_1(d)$ is the same for the two R III moves.  The two R II moves have the same weight $W_1(d)$ and cancel out together.  The partial singularization of the type 1 determines now the partial singularization of the type 7 for the types $r_a$ and $r_b$.  We can write them as singular 3-braids. The signs of the two R III-moves are always opposite by definition. The partial $d$-singularizations are in the brackets.

type 1: $(\sigma_2t_1^+\sigma_2)\sigma_1^{-1}-(\sigma_1t_2^+\sigma_1)\sigma_1^{-1}$

type 7: $\sigma_2(t_1^+\sigma_2\sigma_1^{-1})-\sigma_2(\sigma_2^{-1}\sigma_1t_2^+)$

Hence, the partial $d$-singularization of the type 7 is defined by $t_1^+\sigma_2\sigma_1^{-1}-\sigma_2^{-1}\sigma_1t_2^+$.
\vspace{0,4cm}

But we see that for the global type $r_a$ the weight $W_1(ml)$ changes for the edge $r1-7$,  because the crossing $hm$ does not contribute to this weight.  We compensate the change by adding $(w(hm)-1)/2$ to $W_1(ml)$.  Now the weights are the same. The $ml$-singularizations for the types 1 and 7 are regularly isotopic and cancel out together.
\vspace{0,4cm}

\begin{figure}
\centering
\includegraphics{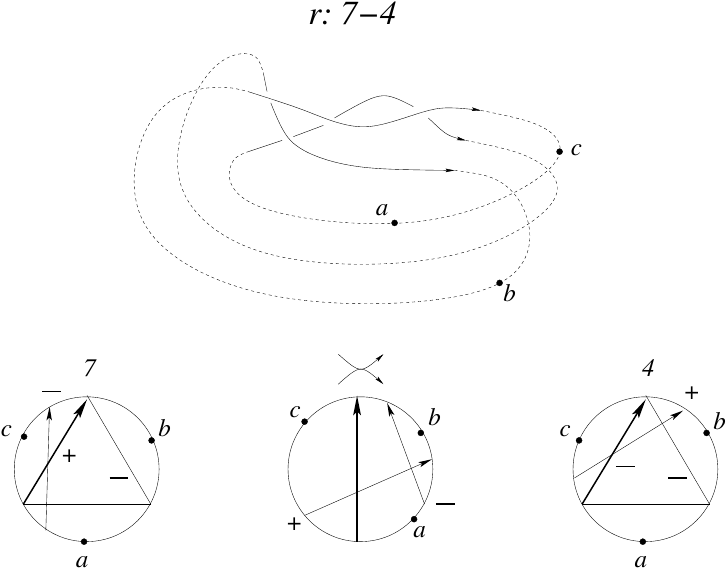}
\caption{\label{r7-4} the edge $r7-4$ in the cube $\Gamma$}  
\end{figure}

The edge $r7-4$ is shown in Fig.~\ref{r7-4}.

We consider the global types types $r_a$ and $r_b$ for the edge $7-4$: Here we have to take care of the signs.  The crossings in the R II moves are just the crossings $d$ in the R III moves and the weights $W_1(d)$ are the same for all R III and R II moves.

R II-move before the R III-move of type 7: $\sigma_2^{-1}\sigma_1(t_2^+)\sigma_2^{-1}-\sigma_2^{-1}\sigma_1\sigma_2(t_2^-)$

type 7: $(t_1^+\sigma_2\sigma_1^{-1})\sigma_2^{-1}-(\sigma_2^{-1}\sigma_1t_2^+)\sigma_2^{-1}$

type 4: $\sigma_1(\sigma_2\sigma_1^{-1}t_2^-)-\sigma_1(t_1^-\sigma_2^{-1}\sigma_1)$

R II-move after the R III-move of type 4: $\sigma_1(t_1^-)\sigma_2^{-1}\sigma_1-(t_1^+)\sigma_1^{-1}\sigma_2^{-1}\sigma_1$

Using only the braid relations we see that the contributions cancel out on the meridian. Hence, the partial $d$-singularization of the type 4 is defined by $t_1^-\sigma_2^{-1}\sigma_1-\sigma_2\sigma_1^{-1}t_2^-$.

The $ml$-singularizations have the same weight, they are regularly isotopic and cancel out together.
\vspace{0,4cm}

\begin{figure}
\centering
\includegraphics{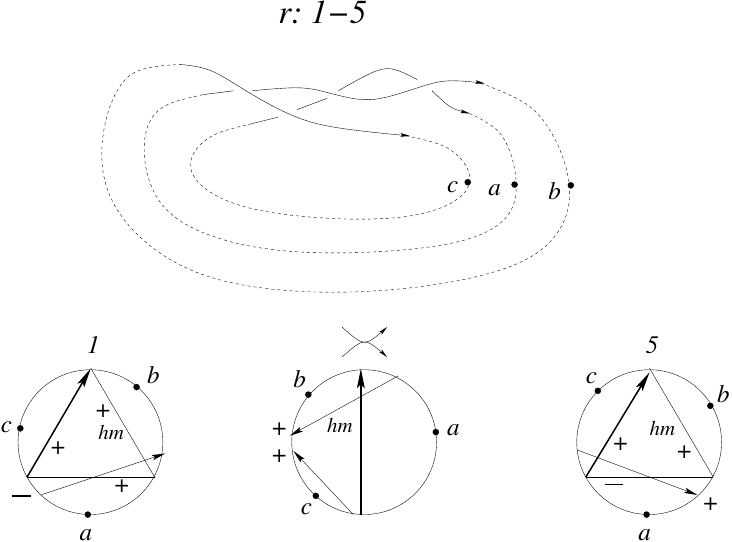}
\caption{\label{r1-5} the edge $r1-5$ in the cube $\Gamma$}  
\end{figure}

The edge $r1-5$ is shown in Fig.~\ref{r1-5}.

The weights $W_1(d)$ are the same for the two triple crossings.

type 1: $(\sigma_2t_1^+\sigma_2)\sigma_2^{-1}-(\sigma_1t_2^+\sigma_1)\sigma_2^{-1}$

This implies the $d$-singularization for the type 5:

type 5: $\sigma_1(\sigma_1^{-1}\sigma_2t_1^+)-\sigma_1(t_2^+\sigma_1\sigma_2^{-1})$
\vspace{0,4cm}

\begin{figure}
\centering
\includegraphics{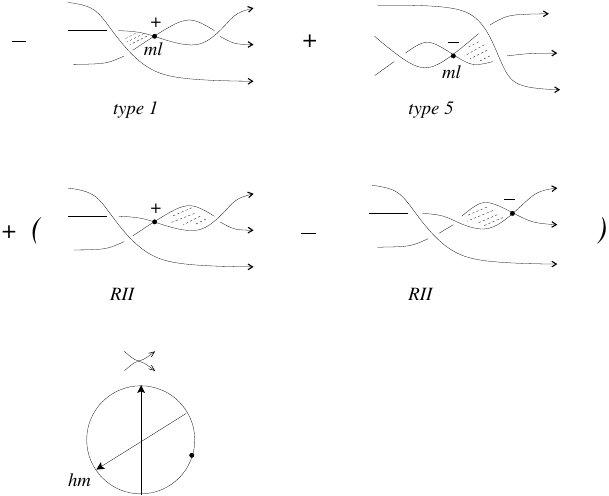}
\caption{\label{ml15} the edge $r1-5$ with the $ml$-singularizations and $d$-singularizations for the R II moves}  
\end{figure}

The sign of the crossing $ml$ changes for the edge $r1-5$.
They have to be singularized for the global type $r_a$, but they are now the same as the crossings $d$ from the R II moves. There is just one branch moving over the auto-tangency.   We see in Fig.~\ref{ml15} that the two tangles with a positive double point are the same and that the two tangles with a negative double point are regularly isotopic.  

However,  the crossing $hm$ of the triple crossings is a f-crossing for the auto-tangency before the R III move of local type 1 (as shown in the figure too, because the triple crossing is of global type $r_a$), but it is no longer a f-crossing for the auto-tangency after the R III move of local type 5. 

The positive double point, as well as the negative double point, enters with the coefficient 

$+x^{W_1(ml)+1} -(x^{W_1(ml)+1}-x^{W_1(ml)})-x^{W_1(ml)}=0$ and hence all contributions cancel out on the meridian of the edge.

This finishes the proof that $\mathbb{L}R_{reg}(x)$ is a 1-cocycle.
\vspace{0,4cm}

The proof carries over for the 1-cocycle for the 2-cable, because we close it abstractly to a circle and we consider as $\infty$ only the starting point of the red knot and we do not care about the endpoint of the red knot, which is the same as the starting point of the black longitude.  In the proof for $\mathbb{L}R_{reg}(x)$ we observe that all cancellations of double points come always from the same crossings in the diagrams. No crossing can move over infinity.  Consequently, we can split the 1-cocycle into two 1-cocycles: singularizations of only crossings of the red knot with itself as well as singularizations of only crossings of the red knot with the black longitude.  Only $0$-crossings can become singular and hence the under-cross has to be black and the over-cross has to be red in this case.

This finishes the proof of the first part of Theorem 1.

\section{The refined tangle equations}

It turns out that it is essential to study  the 1-cocycles  applied to an arc, called the {\em push arc}. 
Let $2K$ be the parallel 2-cable of an auxiliary  framed long knot diagram $K$ and let $2D$ be the parallel 2-cable of a  framed long knot diagram $D$, which we want to study. We denote by $2K2D$ their usual product, i.e.  we glue $2K$ on the left to $2D$. Let $\gamma$ be a regular isotopy of $D$ to another long knot diagram $D'$. We denote the induced regular isotopy of the 2-cables $2K2D$ to $2K2D'$ by $2K\gamma $.

The arc $push(2K,2D)$ starts in $2K2D$ in $M^{red}$. We make $2K$ small and push  it monotonously through $2D$ up to $2D2K$.  Here $D$ is  well defined up to regular isotopy, but $K$ is a fixed diagram.  {\em Attention}: we can not allow to change both $D$ and $K$ by a regular isotopy at the same time, because we evaluate the 1-cocycles only on an arc and not on a loop.
\vspace{0,4cm}

We glue the following arcs to a loop in $M^{red}$ which depends only on $2K$ and $\gamma$:
\vspace{0,4cm}

$push(2K,2D)$  $ (\gamma 2K)$  $(-push(2K,2D'))$  $(- 2K \gamma )$.
\vspace{0,4cm}

The important point is, that this loop is contractible in $M^{red}$. Indeed, $push(2K,2D)$ and $(2K \gamma )$ commute and hence we have a 1-parameter family of loops by shortening $(2K \gamma )$ with respect to the parameter $t$ of the arc $\gamma$. We end up with a loop which is just $push(2K,2D)$  followed by $-push(2K,2D)$ which we shorten  in itself to the point $2K2D$.   

It follows that our 1-cocycle $\mathbb{L}R^{red}_{reg}(x)$  vanishes on the this loop and consequently
\vspace{0,4cm}

$\mathbb{L}R^{red}_{reg}(x)(push(2K,2D))-\mathbb{L}R^{red}_{reg}(x)(push(2K,2D'))=$
\vspace{0,2cm} 

$\mathbb{L}R^{red}_{reg}(x)(2K \gamma)-\mathbb{L}R^{red}_{reg}(x)(\gamma 2K)$.
\vspace{0,2cm}

{\em We denote the corresponding left-hand side of this equation simply by $\Delta(2K)$.  Here we understand $\Delta(2K)$ always as the  element in

 $H_0(M^{red}_+;\mathbb{Z}[x,x^{-1}]) \bigoplus  H_0(M^{red}_-;\mathbb{Z}[x,x^{-1}])$.}

 These are the quantities which we can calculate,  because they do not depend on $\gamma$. But we do not know anything about $\gamma$, not even if it exists at all. Again, by taking the derivative en $x=1$ and applying the HOMFLYPT polynomial to all tangles from the oriented smoothings in the equation we recover the {\em quantum equations for knots} in \cite{FZ} (where we have pushed very simple tangles $T$ through $2D$) in the case of $\mathbb{L}R_{reg}$ for $n=2$.

Clearly,  in $2K \gamma $ we see exactly the same Reidemeister moves as in $- \gamma 2K$, but with opposite signs.  However, it comes from the relative finite type invariants for regular isotopies $W_1(d)$ and $W_1(ml)$, that the contributions of the corresponding Reidemeister moves enter in $\mathbb{L}R^{red}_{reg}(x)(2K \gamma)$ and $\mathbb{L}R^{red}_{reg}(x)( \gamma 2K)$ with different polynomial coefficients and hence in general they will not cancel out pairwise.
\vspace{0,4cm}

We have to express the differences on the right hand side in an explicit way, by using $\gamma$. There are exactly the same Reidemeister moves in $\gamma$ no matter on which end of it is $2K$ and they have the same global types, because $2K$ induce only the trivial permutation of the endpoints of the 2-string link.  But the weights $W_1(d)$ and $W_1(ml)$ depend on the position of $2K$.  

If the knot $2K$ is near to the the starting point $\infty$ of the red oriented  long knot, then we have a small knot $K$ on the red component, a small knot $K$ on the black component and two knots $K$ where for each crossing the over-cross is on one component and the under-cross is on the other component.   Clearly, all $1$-crossings 
of the small knot $K$ on the red component are f-crossings for each Reidemeister move in $2K \gamma$.  Moreover, each crossing of $2K$ with the the under-cross on the red component and the over-cross on the black component is a $1$-crossing and it is also a f-crossing for each Reidemeister move in $2K \gamma$.  On the other hand, if the knot $2K$ is near to the endpoint of the red component, then one easily sees that no $1$-crossing at all of $2K$ is a f-crossing for any Reidemeister move in $\gamma 2K$.

It follows that all the weights $W_1(d)$ and $W_1(ml)$ increase by $w_1(K)+w(K)$ for the Reidemeister moves in $2K \gamma$ with respect to the corresponding moves in $\gamma 2K$.
\vspace{0,4cm}

We have already shown that $\mathbb{L}R^{red}_{reg}(x)$ is a 1-cocycle. Consequently,  $\mathbb{L}R^{red}_{reg}(x)(2\gamma)$ is of the form $ \sum_i D_i(\sum_j  a_{i,j}x^{b_{i,j}})$.

Here $a_{i,j}$ and $b_{i,j}$ are  integers  and $D_i$ are regular isotopy classes of oriented singular 2-string links, each with exactly one ordinary double point which comes from a red-red crossing.  
It follows that 
\vspace{0,2cm}

 $\Delta(2K)=2K \sum_i D_i(\sum_j  a_{i,j}(x^{b_{i,j}+w_1(K)+w(K)}-x^{b_{i,j}}))$.  
 \vspace{0,2cm}

This are the refined tangle equations for long knots with their longitude,  and which moreover split with respect to the positive and the negative double points.  All the same is true for the other case $black-red$ too.

Analogue but even easier arguments give the tangle equations for $n=1$ in Proposition 1.

This finishes the proof of Theorem 1.

\begin{remark}
The tangle equations were generalized with only integer coefficients but with weights of degree 2 in \cite{FMVZ} by using the Gauss diagram formula for the Vassiliev invariant $v_2$ from \cite{PV}.  We conjecture that the tangle equations can be further generalized by using  the Gauss diagram formulas from \cite{CP} for all the Vassiliev invariants which are contained in the HOMFLYPT polynomial.

\end{remark}

\section{Adding the longitude breaks the telescoping effect}

First we explain the {\em telescoping effect} for  $\mathbb{L}R_{reg}(x)$, i.e.  cancellations of contributions for knots without adding the longitude and which makes the invariant trivial.  We show then, that after adding the longitude and considering now $\mathbb{L}R^{red}_{reg}(x)$ the telescoping effect disappears.  This makes $\mathbb{L}R^{red}_{reg}(x)$ potentially to a very powerful tool in order to distinguish knots.  
\vspace{0,4cm}

We push a piece of the auxiliary knot diagram $K$, namely the two curves, over a $0$-crossing of $D$.  This is represented by moving the vertical line under the curves. We are only interested in the double point singularities on the arc from the over-cross of the first crossing to the over-cross of the third crossing of the two curves.  All double points come from $d$-or $ml$-singularizations and therefore all the double points on this arc can come only from the moves involving these three crossings on the curves (double points before or after this arc are blocked by the first and the third crossing and can not move with a knot isotopy into this arc).  The crossings in the curves are always the crossing $hm$ in the R III moves from the moving vertical strand. We have to distinguish two cases. 

\begin{figure}
\centering
\includegraphics{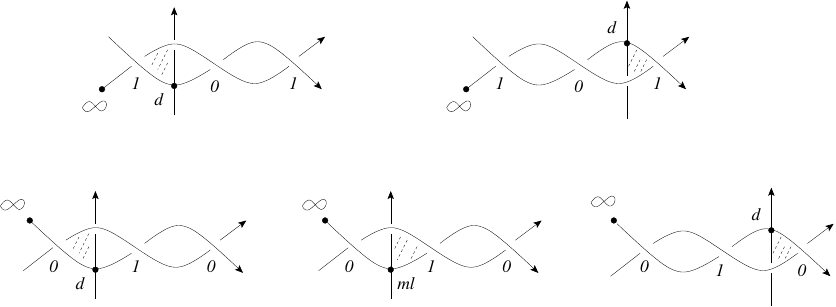}
\caption{\label{tel} the telescoping effect for $n=1$ and moving $K$ over a $0$-crossing of $D$}  
\end{figure}

In the first case, the first crossing in the curves (and hence the third crossing to) is a $1$-crossing and consequently the crossing in the middle is a $0$-crossing.  This is the first line in Fig.~\ref{tel}. 

 It follows that the three R III-moves are of the global types: $r_a,l_b,r_a$.  The move $l_b$ contributes only with $d$-singularizations, but they are not on the arc in which we are interested.  All three moves are of local type $1$ and of positive sign. Let $W_1(d)$ be the weight for the first move $r_a$.  Evidently, the third crossing of the curves, which is a $1$-crossing too,  contributes to $W_1(d)$ of the first move, because when going from $\infty$ to the over-cross of $d$, we pass first through its under-cross. It follows that the two crossings $d$ in the figure share the same $W_1(d)$.  The two tangles are regularly isotopic and they enter $\mathbb{L}R_{reg}(x)$ with different signs (because the first one is on the positive side of the move and the second is on the negative side).  Consequently, the coefficient is $+x^{W_1(d)}-x^{W_1(d)}=0$ and there is no contribution with a double point on this arc.  This is the telescoping effect.

In the second case, the first crossing in the curves (and hence the third crossing to) is a $0$-crossing and consequently the crossing in the middle is a $1$-crossing.  This is the second line in Fig.~\ref{tel}.  The global types are now: $l_b,r_a,l_b$.  Hence the move in the middle contributes with a $ml$-singularization and which is also on the arc we are interested in.  Let $W_1(d)$ be the weight of the first move. Then $W_1(ml)$ in the second move is equal to $W_1(d)$ because the $1$-crossing in the middle does not enter $W_1(ml)$.  For the second move $l_b$ we have now $W_1(d)+1$, because the $1$-crossing in the middle is now a f-crossing.  Again, all three tangles are regularly isotopic and they enter $\mathbb{L}R_{reg}(x)$ with the coefficient from the three moves
\vspace{0,4cm}

$+x^{W_1(d)}+(x^{W_1(d)+1}-x^{W_1(d)})-x^{W_1(d)+1}=0$
\vspace{0,4cm}

This is again the telescoping effect. 

\begin{figure}
\centering
\includegraphics{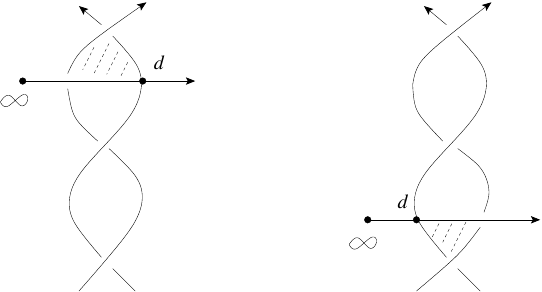}
\caption{\label{second} the telescoping effect for moving $K$ under a $0$-crossing of $D$ }  
\end{figure}

Notice, that when we push $K$ under a $0$-crossing of $D$ then only the global types $r_b$ and $l_b$ occur and hence there is never a $ml$-singularization.  The telescoping effect is even more evident and is shown in Fig.~\ref{second}.  Indeed, moving from $\infty$ to the over-cross of $d$ or $ml$, which is necessarily on a branch of $D$,  we go always only over branches of $K$ and hence no crossing of $K$ can be a f-crossing.

Finally, notice that moving $K$ over a $1$-crossing of $D$ does not contribute at all to 
$\mathbb{L}R_{reg}(x)$. Moving it under a $1$-crossing of $D$ contributes only with a $ml$-singularization for the global type $l_c$.  Hence the double point is just a singularity in the auxiliary knot diagram $K$ and $D$ stays unchanged.

The telescoping effect implies that Proposition 1 is not useful in order to distinguish knots.  However, already $\mathbb{L}R_{reg}$ is useful in order to construct a regular isotopy $\gamma$, compare \cite{FMVZ}.
\vspace{0,4cm}

{\em The crucial point is that the telescoping effect disappears for $\mathbb{L}R^{red}_{reg}(x)$ when we add the black longitude to the red knot. }

We consider the same situation as before and we are only interested in the red-red double points on the thickened red component,  hence we can forget the black vertical strand which is lower  than all other strands and which can be moved away. There are again two cases to consider.

\begin{figure}
\centering
\includegraphics{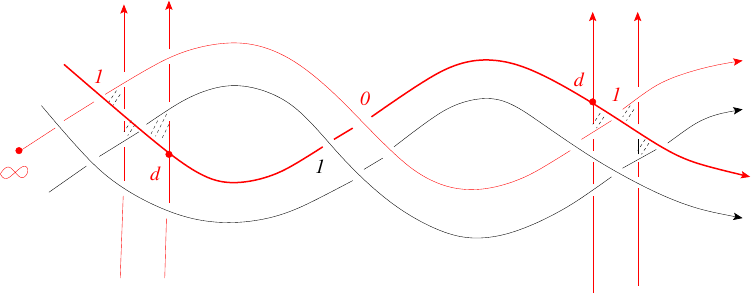}
\caption{\label{nontel1} first case for $n=2$: $r_b,l_b,r_a$}  
\end{figure}

We show the first case in Fig.~\ref{nontel1}.  There are exactly six R III moves for which the crossing $d$ is on the thick red arc from the over-cross of the red $1$-crossing to the over-cross of the other red $1$-crossing. One easily sees that the contributions in the arc from first two R III moves cancel out together, as well as the last two. The third R III move gives a red $d$ as shown in the figure. The middle branch in the move is black and hence the move is necessarily of the global type $r_b$. The red $0$-crossing in the middle gives rise to a move of type $l_b$ and does not contribute with a $ml$-singularization. The fourth move is necessarily of type $r_a$.  Let $W_1(d)$ be the weight of the third move. Again, the $1$-crossing in the fourth move is already a f-crossing for the third move (compare the case $n=1$).  But now there is the foot of a $1$-crossing and hence a f-crossing in the thickened red arc.  Indeed, each crossing with red under-cross and black over-cross is automatically a $1$-crossing.  It follows that the coefficient of the tangle with the double point $d$ on the red thickened arc  is now $+x^{W_1(d)}-x^{W_1(d)+1}$ and it can not cancel out in $\mathbb{L}R^{red}_{reg}(x)$.

\begin{figure}
\centering
\includegraphics{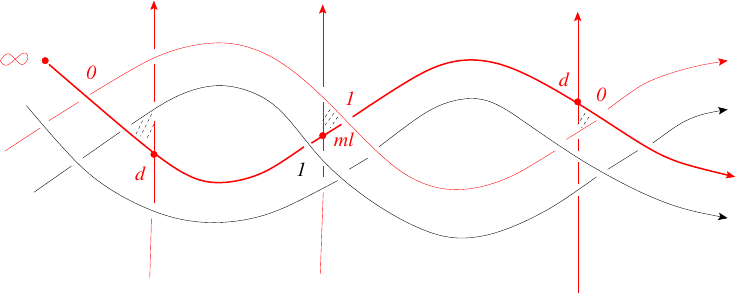}
\caption{\label{nontel2} second case for $n=2$: $r_b,r_a,l_b$}  
\end{figure}

We show the second case in Fig.~\ref{nontel2}.  Besides the previous two $d$-singularizations (the first necessarily of type $r_b$ and the second of type $l_b$) there is now again a $ml$-singularization from the $1$-crossing in the middle, which gives rise to a move of type $r_a$.  Let $W_1(d)$ be the weight of the third move.  But the $1$-crossing with red under-cross and black over-cross on the thickened red arc contributes as a f-crossing to $W_1(ml)$ and hence $W_1(ml)=W_1(d)+1$.  We see in the figure that there are now two $1$-crossings on the arc from the first $d$ to the second $d$ and hence for the second $d$ we have the weight $W_1(d)+2$. All three tangles are regularly isotopic and they enter now with the coefficient

\vspace{0,4cm}

$+x^{W_1(d)}+(x^{W_1(d)+2}-x^{W_1(d)+1})-x^{W_1(d)+2}=+x^{W_1(d)}-x^{W_1(d)+1}$
\vspace{0,4cm}

Notice, that this is exactly the same coefficient as in the first case.  There is still the telescoping effect for moving $2K$ under a red $0$-crossing of $2D$, because the crossings of $2K$ never contribute to the weights here,  compare Fig.~\ref{second}.  Consequently, the singular string link considered above enters $\mathbb{L}R^{red}_{reg}(x)$ always with the coefficient $+x^{W_1(d)}-x^{W_1(d)+1}$ and this coefficient is independent of $w_1(K)+w(K)$.  So, it could well be that for $w_1(K)+w(K)=0$ the above singular string link stays in the invariant $\mathbb{L}R^{red}_{reg}(x)(push(2K,2D))$ if there does not appear a regularly isotopic singular string link with the opposite coefficient in a completely different place in the diagrams $2K 2D$.

All the same is true in the case $\mathbb{L}R^{black-red}_{reg}(x)$ as well.

This suggests that $\mathbb{L}R^{red}_{reg}(x)(push(2K,2D))$ and $\mathbb{L}R^{black-red}_{reg}(x)(push(2K,2D))$ are perhaps  non-trivial singular string link-valued invariants for framed knots if $w_1(K)+w(K)=0$ and we could apply the morphisms $os$ and $cc$ and the machinery of the iterations of our constructions to it. But to confirm this we need again a computer program in order to calculate examples.

Institut de Math\'ematiques de Toulouse, UMR 5219

Universit\'e de Toulouse

31062 Toulouse Cedex 09, France

thomas.fiedler@math.univ-toulouse.fr
\end{document}